\documentclass[a4paper]{amsart}
\usepackage[latin1]{inputenc}
\usepackage{amssymb,amsmath,amsthm}
\usepackage{amsfonts}
\usepackage{ifthen}
\usepackage{paralist}

\theoremstyle{plain}
\newtheorem{theorem}{Theorem}[section]
\newtheorem{proposition}[theorem]{Proposition}
\newtheorem{lemma}[theorem]{Lemma}

\newtheorem{corollary}[theorem]{Corollary}

\newcommand{\vneg}[1]{\ensuremath{#1}^\mathrm{in}} 
\newcommand{\compl}[1]{\ensuremath\overline{#1}} 
\newcommand{\subf}[1][]{\ifthenelse{\equal{#1}{}}{\ensuremath\leq}{\ensuremath\leq_{#1}}}
\newcommand{\psubf}[1][]{\ifthenelse{\equal{#1}{}}{\ensuremath<}{\ensuremath<_{#1}}}
\newcommand{\supf}[1][]{\ifthenelse{\equal{#1}{}}{\ensuremath\geq}{\ensuremath\seq_{#1}}}
\newcommand{\psupf}[1][]{\ifthenelse{\equal{#1}{}}{\ensuremath>}{\ensuremath>_{#1}}}
\newcommand{\incomp}[1][]{\ifthenelse{\equal{#1}{}}{\ensuremath\parallel}{\ensuremath\parallel_{#1}}}
\newcommand{\subc}[1][]{\ifthenelse{\equal{#1}{}}{\ensuremath\preccurlyeq}{\ensuremath\preccurlyeq_{#1}}}
\newcommand{\fequiv}[1][]{\ifthenelse{\equal{#1}{}}{\ensuremath\equiv}{\ensuremath\equiv_{#1}}}
\newcommand{\eqclass}[2][]{\ifthenelse{\equal{#1}{}}{\ensuremath[{#2}]}{\ensuremath[{#2}]_{#1}}}
\newcommand{\subclass}[1][]{\ifthenelse{\equal{#1}{}}{\ensuremath\preccurlyeq}{\ensuremath\preccurlyeq_{#1}}}
\newcommand{\vect}[1]{\ensuremath\mathbf{#1}}
\newcommand{\constf}[1]{\ensuremath \hat{#1}}
\newcommand{\card}[1]{\ensuremath\lvert{#1}\rvert}
\newcommand{\qset}[1]{\ensuremath (\mathcal{O}_A/{\fequiv[{#1}]},\linebreak[0]{\subc[{#1}]})}
\newcommand{\nset}[1]{\ensuremath \underline{#1}}
\DeclareMathOperator{\ea}{ea} 
\DeclareMathOperator{\Ess}{Ess} 
\DeclareMathOperator{\St}{St} 
\DeclareMathOperator{\Int}{Int} 
\DeclareMathOperator{\KER}{ker} 
\DeclareMathOperator{\range}{Im} 
\DeclareMathOperator{\pr}{pr} 
\begin{document}
\hyphenation{Bool-ean clas-ses pa-ram-e-trized}
\title[Subfunctions defined by clones containing all unary operations]{Subfunction relations defined by the clones containing all unary operations}
\author{Erkko Lehtonen}
\address{Institute of Mathematics \\
Tampere University of Technology \\
P.O.~Box 553 \\
FI-33101 Tampere \\
Finland}
\email{erkko.lehtonen@tut.fi}
\begin{abstract}
\noindent For a class $\mathcal{C}$ of operations on a nonempty base set $A$, an operation $f$ is called a $\mathcal{C}$-subfunction of an operation $g$, if $f = g(h_1, \ldots, h_n)$, where all the inner functions $h_i$ are members of $\mathcal{C}$. Two operations are $\mathcal{C}$-equivalent if they are $\mathcal{C}$-subfunctions of each other. The $\mathcal{C}$-subfunction relation is a quasiorder if and only if the defining class $\mathcal{C}$ is a clone. The $\mathcal{C}$-subfunction relations defined by clones that contain all unary operations on a finite base set are examined. For each such clone it is determined whether the corresponding partial order satisfies the descending chain condition and whether it contains infinite antichains.
\end{abstract}
\maketitle

\section{Introduction}
We consider finitary operations on a fixed nonempty base set $A$. Denote by $\mathcal{O}_A$ the set of all operations on $A$. For a class $\mathcal{C} \subseteq \mathcal{O}_A$, we say that $f$ is a \emph{$\mathcal{C}$-subfunction} of $g$, denoted $f \subf[\mathcal{C}] g$, if $f = g(h_1, \ldots, h_n)$, where $h_1, \ldots, h_n \in \mathcal{C}$. Operations $f$ and $g$ are \emph{$\mathcal{C}$-equivalent,} denoted $f \fequiv[\mathcal{C}] g$, if they are $\mathcal{C}$-subfunctions of each other. The relation $\subf[\mathcal{C}]$ is a quasiorder (a reflexive and transitive relation) on $\mathcal{O}_A$ if and only if $\mathcal{C}$ is a clone, i.e., it contains all projections and is closed under functional composition. If $\mathcal{C}$ is a clone, then $\fequiv[\mathcal{C}]$ is indeed an equivalence relation, and $\subf[\mathcal{C}]$ induces a partial order $\subc[\mathcal{C}]$ on the quotient $\mathcal{O}_A / {\fequiv[\mathcal{C}]}$.

$\mathcal{C}$-subfunctions generalize the notion of taking \emph{minors} of functions on finite base sets, an important particular case being that of Boolean functions when the base set is $\{0,1\}$ (see \cite{FH,Pippenger,Wang,WW,Zverovich} for variants). Such generalizations have appeared in many areas of mathematics. For example, $\mathcal{O}_A$-subfunctions were studied by Henno \cite{Henno1971,Henno1977} in the context of Green's equivalences and quasiorders on Menger systems. Equivalences of Boolean functions under actions of the general linear and affine groups of transformations over the two-element field, studied by Harrison \cite{Harrison1964}, correspond to $\mathcal{C}$-equivalences defined by clones of linear Boolean functions, and such linear equivalences have found applications in coding theory and cryptography.

While the lattice of clones on a two-element set was completely described by Post \cite{Post} in the 1940s and his classification has been reproved many times ever since, one of the major ongoing research programs in universal algebra and multi-valued logic is the attempt to understand the structure of the lattice of clones on finite sets with more than two elements. Any approach to Post's theorem is potentially a good candidate for proving Post-like results for large sublattices or sections of the lattice of clones on larger base sets. Representation of classes of functions by forbidden subfunctions playing key role in Zverovich's \cite{Zverovich} proof of Post's theorem, we believe that $\mathcal{C}$-subfunctions could find applications in clone theory and other fields of mathematics.

Previously, we have investigated
the $\mathcal{C}$-subfunction relations defined by the clones of unary, linear and monotone functions on finite base sets in \cite{LinMonot}. We were mostly concerned about the descending chains and antichains, motivated by the fact that representation of classes of functions by minimal sets of forbidden subfunctions is possible if the corresponding $\mathcal{C}$-subfunction partial order satisfies the descending chain condition, and these minimal sets are guaranteed to be finite if the partial order contains only finite antichains. On the other hand, we paid little attention to the closely related ascending chain condition, because it plays no role in the forbidden subfunction characterization.

We now continue our study of $\mathcal{C}$-subfunctions, focusing on the $\mathcal{C}$-subfunction relations defined by the clones that contain all unary operations on a finite set. Such clones constitute a $(k+1)$-element chain $\mathcal{B}_0 \subset \mathcal{B}_1 \subset \dots \subset \mathcal{B}_k$ in the lattice of clones on a set of $k \geq 2$ elements \cite{Burle,Sz}. Our main result is the following.

\begin{theorem}
\label{MainTheorem}
Let $A$ be a finite set with $k \geq 2$ elements, and let $\mathcal{B}_0 \subset \mathcal{B}_1 \subset \dots \subset \mathcal{B}_k$ be the chain of the $k+1$ clones containing all unary operations on $A$. For a particular $0 \leq i \leq k$ we have:

\begin{inparaenum}[(i)]
\item $\subc[\mathcal{B}_i]$ contains infinite descending chains if and only if $2 \leq i \leq k-2$.

\item For $k \geq 3$, $\subc[\mathcal{B}_i]$ contains infinite antichains if and only if $i \leq k-2$.

\item For $k = 2$, $\subc[\mathcal{B}_i]$ contains infinite antichains if and only if $i \leq 1$.
\end{inparaenum}
\end{theorem}

This paper is organized as follows. We present our basic definitions and notation and we outline the proof of Theorem \ref{MainTheorem} in Section \ref{SecDN}. The actual proof amounts to going through all clones $\mathcal{B}_i$ and deciding whether $\subc[\mathcal{B}_i]$ satisfies the descending chain condition and whether it contains infinite antichains. This is done in Sections \ref{SecSL}--\ref{SecAC}.

We analyze the subfunction relations defined by $\mathcal{B}_0$, $\mathcal{B}_k$, $\mathcal{B}_{k-1}$ in Section \ref{SecSL}, and we also give a complete description of the posets $\qset{\mathcal{B}_k}$ and $\qset{\mathcal{B}_{k-1}}$ -- in fact, these posets are finite. In Section \ref{SecBp}, we construct an infinite descending chain of $\mathcal{B}_p$-sub\-func\-tions for $2 \leq p \leq k-2$. Unique standard form representations of quasilinear functions and $\mathcal{C}$-decompositions are introduced in Section \ref{SecB1}, and using these representations we then show that $\subc[\mathcal{B}_1]$ satisfies the descending chain condition. We construct infinite antichains of $\subc[\mathcal{B}_{k-2}]$ in Section \ref{SecAC}.

\section{Definitions and notation}
\label{SecDN}

\subsection{General notation and concepts}

We denote vectors by bold face letters and their components by italic letters, e.g., $\vect{a} = (a_1, \ldots, a_n)$. We also denote the $i$th component of a vector $\vect{a}$ by $\vect{a}(i)$, especially when the vector symbols involve subscripts.

For an integer $n \geq 1$, we denote $\nset{n} = \{0, \ldots, n-1\}$.

For a collection $C$ of disjoint sets, a \emph{transversal} is a set containing exactly one member of each of them. A \emph{partial transversal} of $C$ is a subset of a transversal of $C$.

The \emph{characteristic function} of a subset $S \subseteq A$ is the mapping $\chi_S : A \to \nset{2}$ defined as
\[
\chi_S(x) =
\begin{cases}
1, & \text{if $x \in S$,} \\
0, & \text{if $x \notin S$.}
\end{cases}
\]

\subsection{Functions and clones}

Let $A$ be a fixed nonempty base set. A \emph{function on $A$} is a finitary operation on $A$, i.e., a mapping $f : A^n \rightarrow A$ for some positive integer $n$, called the \emph{arity} of $f$. The set of all functions on $A$ is denoted by $\mathcal{O}_A$.

For a fixed arity $n$, and for $1 \leq i \leq n$, the $n$-ary \emph{$i$th projection,} denoted by $x^n_i$, is the function $(a_1, \ldots, a_n) \mapsto a_i$. The $n$-ary \emph{constant function} having value $a \in A$ everywhere is denoted by $\constf{a}^{n}$. Whenever the arity is clear from the context, we may omit the superscripts indicating arity. The \emph{range,} or \emph{image,} of $f$ is the set $\range f = \{f(\vect{a}) : \vect{a} \in A^n\}$. The \emph{kernel} of $f$ is the equivalence relation $\KER f = \{(\vect{a}, \vect{b}) : f(\vect{a}) = f(\vect{b})\}$ on the domain of $f$.

For $1 \leq i \leq n$, we say that the $i$th variable is \emph{essential} in an $n$-ary function $f$, or $f$ \emph{depends} on the $i$th variable, if there are points $\vect{a} = (a_1, \ldots, a_n)$, $\vect{a}' = (a'_1, \ldots, a'_n)$ such that $a_i \neq a'_i$ and $a_j = a'_j$ for all $j \neq i$ and $f(\vect{a}) \neq f(\vect{a}')$. If a variable is not essential in $f$, then it is \emph{inessential} in $f$. The \emph{essential arity} of $f$, denoted $\ea{f}$, is the number of essential variables in $f$. The \emph{set of essential variables} of $f$ is defined as $\Ess f = \{i : \text{the $i$th variable is essential in $f$}\}$.

If $f$ is an $n$-ary function and $g_1, \ldots, g_n$ are all $m$-ary functions, then the \emph{composition of $f$ with $g_1, \ldots, g_n$,} denoted $f(g_1, \ldots, g_n)$ is an $m$-ary function defined by
\[
f(g_1, \ldots, g_m)(\vect{a}) = f(g_1(\vect{a}), \ldots, g_n(\vect{a})).
\]
This is equivalent to the composition $f \circ g$, where the mapping $g : A^m \rightarrow A^n$ is defined as $g(\vect{a}) = (g_1(\vect{a}), \ldots, g_n(\vect{a}))$, which we simply denote by $g = (g_1, \ldots, g_n)$.

A \emph{class} is a subset $\mathcal{C} \subseteq \mathcal{O}_A$. A \emph{clone} on $A$ is a class $\mathcal{C} \subseteq \mathcal{O}_A$ that contains all projections and is closed under functional composition (i.e., if $f, g_1, \ldots, g_n \in \mathcal{C}$, then $f(g_1, \ldots, g_n) \in \mathcal{C}$ whenever the composition is defined). The clones on $A$ constitute an inclusion-ordered lattice, denoted $\mathcal{L}_A$, where the lattice operations are the following: meet is the intersection, join is the smallest clone containing the union. We denote by $\langle \mathcal{C} \rangle$ the clone generated by $\mathcal{C}$. See \cite{Sz} for a general account on clones.

For any class $\mathcal{C}$, we denote by $\mathcal{C}^{(n)}$ the \emph{$n$-ary part} of $\mathcal{C}$, i.e., $\mathcal{C}^{(n)} = \{f \in \mathcal{C} : \text{$f$ is $n$-ary}\}$. We denote by $\mathcal{C}^{(m,n)}$ the set of mappings of the form $(f_1, \ldots, f_m)$, where each component $f_i$ is a member of $\mathcal{C}^{(n)}$, i.e.,
\[
\mathcal{C}^{(m,n)} = \{(f_1, \ldots, f_m) \in (A^m)^{A^n} : f_1, \ldots, f_m \in \mathcal{C}^{(n)}\}.
\]

For any clone $\mathcal{C}$, $\mathcal{C}^{(1)}$ is a transformation monoid on $A$. Assume that $M$ is an arbitrary transformation monoid on $A$. The \emph{stabilizer} of $M$ is the set
\[
\St M = \{f \in \mathcal{O}_A : \text{$f(g_1, \ldots, g_n) \in M$ whenever $g_1, \ldots, g_n \in M$}\}.
\]
It is easy to see that $\St M$ is a clone, and $\langle M \rangle \subseteq \St M$. An interval of the form $\Int M = [\langle M \rangle; \St M]$ in $\mathcal{L}_A$ is called \emph{monoidal.} It is well-known that $\mathcal{C}^{(1)} = M$ if and only if $\mathcal{C} \in \Int M$ (see, e.g., \cite{Sz}). Thus, $\mathcal{L}_A$ is partitioned in intervals $\Int M$, where $M$ ranges over all submonoids of $\mathcal{O}_A^{(1)}$.

It is well-known that for a finite base set $A$ with $k \geq 2$ elements, the monoidal interval $\Int \langle \mathcal{O}_A^{(1)} \rangle$ is the $(k + 1)$-element chain
\[
\langle \mathcal{O}_A^{(1)} \rangle = \mathcal{B}_0 \subset \mathcal{B}_1 \subset \mathcal{B}_2 \subset \dots \subset \mathcal{B}_{k-1} \subset \mathcal{B}_k = \mathcal{O}_A,
\]
where the clones $\mathcal{B}_i$ are defined as follows. For $2 \leq i \leq k$, $\mathcal{B}_i$ consists of all essentially at most unary functions and all functions whose range contains at most $i$ elements. $\mathcal{B}_1$ consists of all essentially at most unary functions and all \emph{quasilinear} functions, i.e., functions having the form $g(h_1(x_1) \oplus \dots \oplus h_n(x_n))$ with $h_1, \ldots, h_n : A \rightarrow \nset{2}$, $g : \nset{2} \rightarrow A$ arbitrary mappings and $\oplus$ denoting addition modulo $2$. The fact that $\mathcal{B}_{k-1}$ is the unique maximal clone that contains all unary functions was first proved by S\l{}upecki \cite{Slupecki}. The unrefinable chain $\mathcal{B}_2 \subset \mathcal{B}_3 \subset \dots \subset \mathcal{B}_k$ was also known to him. The description of the monoidal interval was completed by Burle \cite{Burle} who proved that $\mathcal{B}_1$ is the only other clone properly containing $\mathcal{B}_0$.

\subsection{$\mathcal{C}$-subfunctions}

Let $\mathcal{C}$ be a class of functions on $A$. We say that a function $f$ is a \emph{$\mathcal{C}$-subfunction} of a function $g$, denoted $f \subf[\mathcal{C}] g$, if $f = g(h_1, \ldots, h_m)$ for some $h_1, \ldots, h_m \in \mathcal{C}$, i.e., $f = g \circ h$ for some $h \in \mathcal{C}^{(m,n)}$ where $m$ and $n$ are the arities of $g$ and $f$, respectively. If $f$ and $g$ are $\mathcal{C}$-subfunctions of each other, we say that they are \emph{$\mathcal{C}$-equivalent} and denote $f \fequiv[\mathcal{C}] g$. If $f \subf[\mathcal{C}] g$ but $g \not\subf[\mathcal{C}] f$, we say that $f$ is a \emph{proper $\mathcal{C}$-subfunction} of $g$ and denote $f \psubf[\mathcal{C}] g$. If both $f \not\subf[\mathcal{C}] g$ and $g \not\subf[\mathcal{C}] f$, we say that $f$ and $g$ are \emph{$\mathcal{C}$-incomparable} and denote $f \parallel_\mathcal{C} g$. If the class $\mathcal{C}$ is clear from the context, we may simplyfy the notation and omit the subscripts indicating the class.

We have now defined families of binary relations $\subf[\mathcal{C}]$ and $\fequiv[\mathcal{C}]$ on $\mathcal{O}_A$, parametrized by the class $\mathcal{C}$. The following basic properties of $\subf[\mathcal{C}]$ and $\fequiv[\mathcal{C}]$ can be proved easily.

For any class $\mathcal{C}$, the set of $\mathcal{C}$-subfunctions of $x_1$ is $\mathcal{C}$, and therefore the relations $\subf[\mathcal{C}]$ and $\subf[\mathcal{K}]$ are distinct for $\mathcal{C} \neq \mathcal{K}$. Also, for any classes $\mathcal{C}$ and $\mathcal{K}$, $\subf[\mathcal{C}]$ is a subrelation of $\subf[\mathcal{K}]$ if and only if $\mathcal{C} \subseteq \mathcal{K}$. For any clones $\mathcal{C}$ and $\mathcal{K}$, $\fequiv[\mathcal{C}]$ is a subrelation of $\fequiv[\mathcal{K}]$ whenever $\mathcal{C} \subseteq \mathcal{K}$. However, it is possible that $\fequiv[\mathcal{C}]$ and $\fequiv[\mathcal{K}]$ coincide even if $\mathcal{C} \neq \mathcal{K}$.

The relation $\subf[\mathcal{C}]$ is reflexive if and only if the class $\mathcal{C}$ contains all projections; and $\subf[\mathcal{C}]$ is transitive if and only if $\mathcal{C}$ is closed under functional composition. Hence, $\subf[\mathcal{C}]$ is a quasiorder on $\mathcal{O}_A$ if and only if $\mathcal{C}$ is a clone. If $\mathcal{C}$ is a clone, then $\fequiv[\mathcal{C}]$ is an equivalence relation. The $\fequiv[\mathcal{C}]$-class of $f$ is denoted by $\eqclass[\mathcal{C}]{f}$. As for quasiorders, $\subf[\mathcal{C}]$ induces a partial order $\subc[\mathcal{C}]$ on $\mathcal{O}_A / {\fequiv[\mathcal{C}]}$.

It is clear that $\range f \subseteq \range g$ for any $f \subf[\mathcal{C}] g$ and any $\mathcal{C}$. Therefore, any $\mathcal{C}$-equiv\-a\-lent functions have the same range. This implies in particular that for any element $a \in A$, the constant functions $\constf{a}$ of all arities form a $\fequiv[\mathcal{C}]$-class for any clone $\mathcal{C}$, and these classes are minimal in the partial order $\subc[\mathcal{C}]$ of $\mathcal{O}_A/{\fequiv[\mathcal{C}]}$.

In what follows, we assume that the base set $A$ is finite and $\card{A} = k \geq 2$. Because it is immaterial what the elements of the base set are, we assume, without loss of generality, that $A = \{0, 1, \ldots, k - 1\} = \nset{k}$.

\subsection{Proof of Theorem \ref{MainTheorem}}

In the subsequent sections, we will decide for each clone $\mathcal{B}_i$ ($0 \leq i \leq k$) whether $\subc[\mathcal{B}_i]$ satisfies the descending chain condition and whether it contains infinite antichains. The statement about the descending chain condition follows from Theorems \ref{B0chain}, \ref{BkStr}, \ref{Bk-1Str}, \ref{Bpchain}, \ref{B1chain}. The statement about antichains for $k \geq 3$ follows from Theorems \ref{Bk-1Str}, \ref{Bk-2antichain}, taking into account that ${\subf[\mathcal{C}]} \subseteq {\subf[\mathcal{K}]}$ if and only if $\mathcal{C} \subseteq \mathcal{K}$. For $k = 2$, the statement about antichains follows from Theorem \ref{BkStr} and the fact that $\mathcal{B}_1$ contains infinite antichains, proved in \cite[Theorem 5.10]{LinMonot}.

\section{The smallest and largest clones in $\Int \langle \mathcal{O}^{(1)}_A \rangle$}
\label{SecSL}

\subsection{$\mathcal{B}_0$-subfunctions}

Denote by $\mathcal{J}_A$ the clone of all projections on $A$. It is clear that every nonconstant function $f$ is $\mathcal{J}_A$-equivalent (and hence $\mathcal{C}$-equivalent for any clone $\mathcal{C}$) to the function $f^\mathrm{ess}$ of arity $\ea f$, obtained by deleting all inessential variables of $f$. We can also agree that $(\constf{a}^n)^\mathrm{ess} = \constf{a}^1$.

\begin{lemma}
\label{ea}
Let $M$ be a transformation monoid on $A$. If $f \subf[\langle M \rangle] g$, then $\ea f \leq \ea g$.
\end{lemma}
\begin{proof}
Let $f = g(h_1, \ldots, h_m)$ for some $h_1, \ldots, h_m \in \langle M \rangle$. Each essential variable of $f$ has to be essential in at least one of the inner functions $h_i$ substituted for an essential variable of $g$. Since the $h_i$'s are essentially at most unary, it is clear that $\ea f \leq \ea g$.
\end{proof}

\begin{proposition}
\label{proptm}
For any transformation monoid $M$ on $A$, there is no infinite descending chain of $\langle M \rangle$-subfunctions.
\end{proposition}
\begin{proof}
Suppose, on the contrary, that there is an infinite descending chain
\[
f_1 \psupf[\langle M \rangle] f_2 \psupf[\langle M \rangle] f_3 \psupf[\langle M \rangle] \cdots.
\]
Since each $f_i$ is $\langle M \rangle$-equivalent to $f_i^\mathrm{ess}$, we can assume that all variables are essential in $f_i$. Lemma \ref{ea} implies that there is an $m$ such that all functions $f_i$ with $i \geq m$ have the same arity. We have reached a contradiction, because there are only a finite number of functions of any fixed arity.
\end{proof}

\begin{theorem}
\label{B0chain}
$\subc[\mathcal{B}_0]$ satisfies the descending chain condition.
\end{theorem}
\begin{proof}
A special case of Proposition \ref{proptm}.
\end{proof}

\subsection{The structure of $\subc[\mathcal{B}_k]$}

For a class $\mathcal{C}$ and a subset $S \subseteq A$, we denote $\mathcal{C}_S = \{f \in \mathcal{C} : \range f = S\}$.

\begin{proposition}
\label{proposition:AllUnary}
Suppose that $\mathcal{C} \in \Int \langle \mathcal{O}_A^{(1)} \rangle$. Then for every nonempty subset $S$ of $A$, every nonempty $\mathcal{C}_S$ is a $\fequiv[\mathcal{C}]$-class.
\end{proposition}
\begin{proof}
Let $\xi$ be any fixed unary function in $\mathcal{C}_S$ such that $\xi(a) = a$ for all $a \in S$, and let $f \in \mathcal{C}_S$ be $n$-ary. It is clear that $\xi(f) = f$, so $f \subf[\mathcal{C}] \xi$. For each $b \in S$, let $\vect{u}_b \in f^{-1}(b)$. For $i = 1, \ldots, n$, define the unary function $g_i$ as $g_i(a) = (\vect{u}_{\xi(a)})_i$. Then
\[
f(g_1, \ldots, g_n)(a) = f(g_1(a), \ldots, g_n(a)) = f(\vect{u}_{\xi(a)}) = \xi(a),
\]
so $\xi = f(g_1, \ldots, g_n)$. Since $\mathcal{C}$ contains all unary functions, $\xi \subf[\mathcal{C}] f$. We have shown that all functions in $\mathcal{C}_S$ are $\mathcal{C}$-equivalent to $\xi$. By the transitivity of $\fequiv[\mathcal{C}]$, the members of $\mathcal{C}_S$ are pairwise $\mathcal{C}$-equivalent. The fact that $\mathcal{C}$-equivalent functions have the same range now implies that $\mathcal{C}_S$ is a $\fequiv[\mathcal{C}]$-class.
\end{proof}

Proposition \ref{proposition:AllUnary} gives a complete characterization of the $\fequiv[\mathcal{B}_k]$-classes: $f \fequiv[\mathcal{B}_k] g$ if and only if $\range f = \range g$. Since $A$ is finite, there are only a finite number of equivalence classes, and therefore there simply cannot exist an infinite descending chain of $\mathcal{B}_k$-subfunctions nor an infinite antichain of $\mathcal{B}_k$-incomparable functions. In fact, it is easy to see that $\qset{\mathcal{B}_k}$ is isomorphic to $(\mathcal{P}(A) \setminus \{\emptyset\}, {\subseteq})$, the power set lattice of $A$ with the bottom element removed; this was in fact proved by Henno \cite{Henno1971} who studied Green's relations on Menger systems. The largest chain of this poset has $k$ elements, and by Sperner's theorem \cite{Sperner}, the largest antichain has $\binom{k}{\lfloor k/2 \rfloor}$ elements.

\begin{theorem}
\label{BkStr}
The poset $\qset{\mathcal{B}_k}$ is isomorphic to $(\mathcal{P}(A) \setminus \{\emptyset\}, {\subseteq})$. The largest chain of this poset has $k$ elements and the largest antichain has $\binom{k}{\lfloor k/2 \rfloor}$ elements.
\end{theorem}

\subsection{The structure of $\subc[\mathcal{B}_{k-1}]$}

Let $\alpha$, $\beta$, $\gamma$ be distinct elements of $A$. We say that $(\alpha, \beta, \gamma)$ is an \emph{essential triple} for an $n$-ary function $f$ if there exist $\vect{a}, \vect{b}, \vect{c} \in A^n$ and $1 \leq i \leq n$ such that $a_j = b_j$ for every $j \neq i$, $a_i = c_i$, and $f(\vect{a}) = \alpha$, $f(\vect{b}) = \beta$, $f(\vect{c}) = \gamma$.

The following lemma is due to Mal'tsev \cite{Maltsev} who slightly improved earlier results by Yablonski \cite{Jablonski} and Salomaa \cite{Salomaa} (see also \cite{Rosenberg}).

\begin{lemma}
\label{lemma:esstriple}
If $f$ has at least two essential variables and takes on more than two values, then $f$ possesses an essential triple. Conversely, if $f$ possesses an essential triple, then $f$ has at least two essential variables and each value of $f$ appears in an essential triple for $f$.
\end{lemma}

For $1 \leq i \leq n$, we define the \emph{projection} of a subset $B \subseteq A^n$ onto its $i$th component by
\[
\pr_i B = \{b_i : (b_1, \ldots, b_n) \in B\}.
\]

\begin{corollary}
\label{restr}
Let $f$ be an $n$-ary function with at least two essential variables and $\card{\range f} = r \geq 3$. Then there is a transversal $B$ of $\KER f$ such that $\card{\pr_i B} < r$ for $1 \leq i \leq n$.
\end{corollary}
\begin{proof}
By Lemma~\ref{lemma:esstriple}, there are elements $\alpha, \beta, \gamma \in A$ and points $\vect{a}, \vect{b}, \vect{c} \in A^n$ such that the conditions for an essential triple for $f$ are satisfied. Choose any $r-3$ points $\vect{d}_4, \ldots, \vect{d}_r$ such that $B = \{\vect{a}, \vect{b}, \vect{c}, \vect{d}_4, \ldots, \vect{d}_r\}$ is a transversal of $\KER f$. It is clear that $\card{\pr_i B} \leq r-1$ for $1 \leq i \leq n$.
\end{proof}

Now we can characterize the $\fequiv[\mathcal{B}_{k-1}]$-classes. The classes contained in $\mathcal{B}_{k-1}$ are given by Proposition \ref{proposition:AllUnary}, and we will show that $\mathcal{O}_A \setminus \mathcal{B}_{k-1}$ is a $\fequiv[\mathcal{B}_{k-1}]$-class. For this purpose, we only have to show that $f \fequiv[\mathcal{B}_{k-1}] g$ for any $f, g \in \mathcal{O}_A \setminus \mathcal{B}_{k-1}$.

\begin{proposition}
If $f, g \in \mathcal{O}_A \setminus \mathcal{B}_{k-1}$, then $f \fequiv[\mathcal{B}_{k-1}] g$.
\end{proposition}
\begin{proof}
Let $f$ be an $n$-ary function and $g$ an $m$-ary function, both in $\mathcal{O}_A \setminus \mathcal{B}_{k-1}$. Hence $f$ and $g$ are essentially at least binary and $\range f = \range g = A$. By Corollary \ref{restr}, there is a $k$-element set $B = \{\vect{d}_0, \ldots, \vect{d}_{k-1}\} \subseteq A^m$ such that $g(\vect{d}_a) = a$ for every $a \in A$ and $\card{\pr_i B} \leq k-1$ for $1 \leq i \leq n$. Let $h : A^n \to A^m$ be defined as $h(\vect{a}) = \vect{d}_{f(\vect{a})}$. We clearly have that $f = g \circ h$ and $h \in \mathcal{B}_{k-1}^{(m,n)}$. Thus $f \subf[\mathcal{B}_{k-1}] g$. Similarly, we can show that $g \subf[\mathcal{B}_{k-1}] f$.
\end{proof}

There are only a finite number of $\fequiv[\mathcal{B}_{k-1}]$-classes, and therefore $\subc[\mathcal{B}_{k-1}]$ contains only finite chains and antichains. The structure of the poset $\qset{\mathcal{B}_{k-1}}$ can easily be described in more detail. We denote by $P \oplus Q$ the \emph{linear sum} of posets $P$ and $Q$, and we denote by $\mathbf{1}$ the one-element chain. (See \cite{DP} for more details.) Then $\qset{\mathcal{B}_{k-1}}$ is isomorphic to $(\mathcal{P}(A) \setminus \{\emptyset\}, {\subseteq}) \oplus \mathbf{1}$, the power set lattice of $A$ with the bottom element removed and a new top element added.

\begin{theorem}
\label{Bk-1Str}
The poset $\qset{\mathcal{B}_{k-1}}$ is isomorphic to $(\mathcal{P}(A) \setminus \{\emptyset\}, {\subseteq}) \oplus \mathbf{1}$. The largest chain of this poset has $k+1$ elements and the largest antichain has $\binom{k}{\lfloor k/2 \rfloor}$ elements.
\end{theorem}

\section{Descending chains in $\subc[\mathcal{B}_p]$ for $2 \leq p \leq k-2$}
\label{SecBp}

In this section, we assume that $\card{A} = k \geq 4$. For $n \geq 3$, define the $(n+1)$-ary function $f_n$ on $A$ as
\[
f_n(\vect{a}) =
\begin{cases}
a, & \text{if $a_1 = \dots = a_n = a$ and $a_{n+1} = 0$,} \\
1, & \text{if $\vect{a} \in \{\vect{u}_n, \vect{v}_n, \vect{w}_n\}$,} \\
0, & \text{otherwise,}
\end{cases}
\]
where the $(n+1)$-vectors $\vect{u}_n$, $\vect{v}_n$, $\vect{w}_n$ are defined recursively as
\begin{align*}
&\vect{u}_3 = (1,1,2,0), & &\vect{v}_3 = (2,1,2,0), & &\vect{w}_3 = (3,2,1,0), \\
&\vect{u}_4 = (1,2,1,2,0), & &\vect{v}_4 = (2,3,2,1,0), & &\vect{w}_4 = (3,2,1,2,0), \\
\intertext{and for $n \geq 5$ and $2 \leq i \leq n+1$,}
&\vect{u}_n(1) = 1 & &\vect{v}_n(1) = 2, & &\vect{w}_n(1) = 3, \\
&\vect{u}_n(i) = \vect{v}_{n-1}(i-1), & &\vect{v}_n(i) = \vect{u}_{n-1}(i-1), & &\vect{w}_n(i) = \vect{u}_{n-1}(i-1).
\end{align*}

\begin{proposition}
For $2 \leq p \leq k-2$ and for any $n \geq 3$, $f_{n+1} \psubf[\mathcal{B}_p] f_n$.
\end{proposition}
\begin{proof}
We first observe that $f_{n+1} = f_n(x_2, x_3, \ldots, x_n, x_{n+1}, g)$, where the $(n+2)$-ary function $g$ is defined as
\[
g(\vect{a}) =
\begin{cases}
0, & \text{if $\vect{a} \in \{\vect{u}_{n+1}, \vect{v}_{n+1}, \vect{w}_{n+1}\}$ or $a_1 = \dots = a_{n+1}$,} \\
1, & \text{otherwise.}
\end{cases}
\]
Since all projections and $g$ are members of $\mathcal{B}_2$, we conclude that $f_{n+1} \subf[\mathcal{C}] f_n$ for every $\mathcal{C} \supseteq \mathcal{B}_2$.

We then show that $f_n \not\subf[\mathcal{B}_{k-2}] f_{n+1}$ and hence $f_n \not\subf[\mathcal{C}] f_{n+1}$ for every $\mathcal{C} \subseteq \mathcal{B}_{k-2}$. Suppose, on the contrary, that $f_n \subf[\mathcal{B}_{k-2}] f_{n+1}$. Then there exist $(n+1)$-ary functions $h_1, \ldots, h_{n+2} \in \mathcal{B}_{k-2}$ such that $f_n = f_{n+1}(h_1, \ldots, h_{n+2})$. Let us denote $h = (h_1, \ldots, h_{n+2})$, and for $a \in A$ and $n \geq 1$, denote by $\vect{e}^n_a$ the $(n+1)$-vector whose first $n$ components are equal to $a$ and the last component is $0$. We clearly have that for $a \in \{2, \ldots, k-1\}$, $h(\vect{e}^n_a) = \vect{e}^{n+1}_a$ and $h[\{\vect{e}^n_1, \vect{u}_n, \vect{v}_n, \vect{w}_n\}] \subseteq \{\vect{e}^{n+1}_1, \vect{u}_{n+1}, \vect{v}_{n+1}, \vect{w}_{n+1}\}$.

We say that a function $f$ is a \emph{projection beyond $a \in A$} if $\Ess f = \{i\}$ for some $i$ and $f(a_1, \ldots, a_n) = a_i$ whenever $a_i \neq a$. Note that all projections are projections beyond every $a \in A$.

If $h(\vect{e}^n_1) = \vect{u}_{n+1}$, then $\card{\range h_i} \geq k-1$ whenever $\vect{u}_{n+1}(i) = 1$ and so $h_i$ is essentially unary and hence a projection beyond $0$. Apart from $h_{n+2}$, each of the other inner functions $h_i$ is either essentially at least binary, and therefore $0,1 \notin \range h_i$, or $h_i$ is essentially unary such that $h_i(j, \ldots, j) = j$ for $2 \leq j \leq k-1$ but $h_i(1, \ldots, 1) \neq 1$. Whatever the case may be, it is only possible that $h$ maps both $\vect{u}_n$ and $\vect{v}_n$ to $\vect{u}_{n+1}$. However, this is an impossibility, because $\vect{u}_n$ and $\vect{v}_n$ do not have $1$'s in common positions.

Similarly, we deduce that it is not possible that $h(\vect{e}^n_1) = \vect{v}_{n+1}$ or $h(\vect{e}^n_1) = \vect{w}_{n+1}$. We are only left with the case that $h(\vect{e}^n_1) = \vect{e}^{n+1}_1$. Then for $1 \leq i \leq n+1$, $\card{\range h_i} \geq k-1$ and hence $h_i$ is essentially unary. Moreover, $h_i$ is a projection beyond $0$ and $n+1 \notin \Ess h_i$.

We observe that $\vect{u}_n(1) = 1$, $\vect{v}_n(1) = 2$, $\vect{w}_n(1) = 3$ and for $2 \leq i \leq n$, $\card{\{\vect{u}_n(i), \vect{v}_n(i), \vect{w}_n(i)\}} = 2$. Since $h_1$ is a projection beyond $0$, it is only possible that $\Ess h_1 = \{1\}$, and we can deduce that $h(\vect{v}_n) = \vect{v}_{n+1}$, $h(\vect{w}_n) = \vect{w}_{n+1}$ and either $h(\vect{u}_n) = \vect{u}_{n+1}$ or $h(\vect{u}_n) = \vect{e}^{n+1}_1$. From the recursive definition, we see that at most one component of $\vect{v}_{n+1}$ equals $3$. If $\vect{v}_{n+1}$ contains a $3$, then $\vect{v}_n$ does not, and there is no way we could produce a $3$ from the vector $\vect{v}_n$ that does not contain a $3$ by a projection beyond $0$ with the $(n+1)$-th variable inessential. Otherwise, $\vect{u}_{n+1}$ contains a $3$ but $\vect{u}_n$ does not, and we deduce in a similar way that $h(\vect{u}_n) \neq \vect{u}_{n+1}$, so we have that $h(\vect{u}_n) = \vect{e}^{n+1}_1$. Since $\vect{u}_n$ and $\vect{v}_n$ do not have $1$'s at common positions, $h(\vect{v}_n)$ would be a vector with no $1$'s, a contradiction. This concludes the proof that $f_n \not\subf[\mathcal{B}_{k-2}] f_{n+1}$.
\end{proof}

\begin{theorem}
\label{Bpchain}
For $2 \leq p \leq k-2$, there is an infinite descending chain of $\mathcal{B}_p$-sub\-func\-tions.
\end{theorem}

\section{$\mathcal{B}_1$-subfunctions}
\label{SecB1}

In this section, we will show that $\subc[\mathcal{B}_1]$ satisfies the descending chain condition. To this end, we first introduce two technical notions that we will need in our analysis, namely unique standard form representations of quasilinear functions and $\mathcal{C}$-decompositions.

\subsection{Unique representations of quasilinear functions}

Functions of the form $f = g(h_1(x_1) \oplus \dots \oplus h_n(x_n))$, where $h_1, \ldots, h_n : A \rightarrow \nset{2}$, $g : \nset{2} \rightarrow A$ and $\oplus$ denotes addition modulo $2$, are called \emph{quasilinear.} The mappings $h_i$ are in fact characteristic functions of subsets $S_i \subseteq A$. Then $h_1(x_1) \oplus \dots \oplus h_n(x_n)$ is the characteristic function of the set
\[
S'_1 \bigtriangleup S'_2 \bigtriangleup \dots \bigtriangleup S'_n \subseteq A^n,
\]
where $\bigtriangleup$ denotes symmetric difference and
\[
S'_i = \{(a_1, \ldots, a_n) \in A^n : a_i \in S_i\} = A^{i-1} \times S_i \times A^{n-i} \subseteq A^n.
\]

The \emph{negation} $\compl{h}$ of a mapping $h : A^n \to \nset{2}$ is defined as $\compl{h}(\vect{a}) = h(\vect{a}) \oplus 1$ for all $\vect{a} \in A^n$. The \emph{inner negation} $\vneg{g}$ of a mapping $g : \nset{2} \to A$ is defined as $\vneg{g}(b) = g(b \oplus 1)$ for each $b \in \nset{2}$. We note that the negation of the characteristic function of a subset $S \subseteq A$ is the characteristic function of the complement of $S$.

By the preceding remarks, it is not difficult to see that the representation of a nonconstant quasilinear function in the form $f = g(h_1(x_1) \oplus \dots \oplus h_n(x_n))$ is unique up to the negation of some of the functions $h_i$ and the inner negation of $g$ if the number of negated $h_i$'s is odd.

For any fixed element $a \in A$, we can choose all the functions $h_i$ such that they are characteristic functions of subsets of $A$ that do not contain $a$, and this way we achieve unique representations of nonconstant quasilinear functions. We say that the representation $f = g(h_1(x_1) \oplus \dots \oplus h_n(x_n))$ of a nonconstant quasilinear function $f$ is in \emph{standard form} if $h_i(0) = 0$ for every $i = 1, \ldots, n$. Standard forms are unique.

\subsection{$\mathcal{C}$-decompositions}

Let $\mathcal{C}$ be a clone. If $f = g(\phi_1, \ldots, \phi_m)$ for $\phi_1, \ldots, \phi_m \in \mathcal{C}$, we say that the $(m+1)$-tuple $(g, \phi_1, \ldots, \phi_m)$ is a \emph{$\mathcal{C}$-decomposition} of $f$. We often avoid referring explicitly to the tuple and we simply say that $f = g(\phi_1, \ldots, \phi_m)$ is a $\mathcal{C}$-decomposition. $\mathcal{C}$-decompositions always exist for all clones $\mathcal{C}$ and all functions $f$, because $f = f(x_1, \ldots, x_n)$ and projections are members of every clone. We call a $\mathcal{C}$-decomposition $(g, \phi_1, \ldots, \phi_m)$ of a nonconstant function $f$ \emph{minimal,} if the number $m$ of inner functions is the smallest possible among all $\mathcal{C}$-decompositions of $f$, and we call this smallest number the \emph{$\mathcal{C}$-degree} of $f$, denoted $\deg_\mathcal{C} f$. We agree that the $\mathcal{C}$-degree of a constant function is $0$. It is clear that $\deg_\mathcal{C} f \leq \ea f$ for any function $f$.

\begin{lemma}
\label{Cdegree}
If $f \subf[\mathcal{C}] g$ then $\deg_\mathcal{C} f \leq \deg_\mathcal{C} g$.
\end{lemma}
\begin{proof}
Let $g = s(\phi_1, \ldots, \phi_d)$ be a minimal $\mathcal{C}$-decomposition. Since $f \subf[\mathcal{C}] g$, we have that $f = g(h_1, \ldots, h_m)$ for some $h_1, \ldots, h_m \in \mathcal{C}$. Then
\[
f = s(\phi_1, \ldots, \phi_d)(h_1, \ldots, h_m) = s(\phi'_1, \ldots, \phi'_d),
\]
where $\phi'_i = \phi_i(h_1, \ldots, h_m) \in \mathcal{C}$. Thus, $\deg_\mathcal{C} f \leq d$. The claim also holds for constant functions, because all $\mathcal{C}$-subfunctions of a constant function are constant.
\end{proof}

\begin{corollary}
$\mathcal{C}$-equivalent functions have the same $\mathcal{C}$-degree.
\end{corollary}

An $m$-tuple ($m \geq 2$) $(\phi_1, \ldots, \phi_m)$ of $n$-ary functions is \emph{functionally dependent,} if there is an $(m-1)$-ary function $g$ and an $i$ such that $\phi_i = g(\phi_1, \ldots, \phi_{i-1},\linebreak[0] \phi_{i+1},\linebreak[0] \ldots,\linebreak[0] \phi_m)$. A tuple $(\phi_1, \ldots, \phi_m)$ is \emph{functionally independent} if it is not functionally dependent. We often omit the tuple notation and we simply say that functions $\phi_1, \ldots, \phi_m$ are functionally dependent or independent.

\begin{lemma}
In a minimal $\mathcal{C}$-decomposition $(g, \phi_1, \ldots, \phi_d)$ of $f$, the inner functions $\phi_1, \ldots, \phi_d$ are functionally independent.
\end{lemma}
\begin{proof}
Suppose, on the contrary, that there is a $(d-1)$-ary function $h$ and an $i$ such that $\phi_i = h(\phi_1, \ldots, \phi_{i-1}, \phi_{i+1}, \ldots, \phi_d)$. Then
\[
f = g(\phi_1, \ldots, \phi_d)
= g(x_1, \ldots, x_{i-1}, h, x_i, \ldots, x_{d-1})(\phi_1, \ldots, \phi_{i-1}, \phi_{i+1}, \ldots, \phi_d),
\]
a contradiction to the minimality of $(g, \phi_1, \ldots, \phi_d)$.
\end{proof}

Any $m$-tuple ($m \geq 2$) of functions containing a constant function is clearly functionally dependent, and therefore none of the inner functions of a minimal $\mathcal{C}$-decomposition is a constant function. The following more general statement also holds.

\begin{lemma}
\label{Cdecev}
If $f = s(\phi_1, \ldots, \phi_d)$ is a minimal $\mathcal{C}$-decomposition, then for every $m \geq 1$ and for all $m$-element subsets $S \subseteq \{1, \ldots, d\}$,
\[
\card{\bigcup_{i \in S} \Ess \phi_i} \geq m.
\]
\end{lemma}
\begin{proof}
For the sake of contradiction, assume without loss of generality that
\[
\bigcup_{i=1}^m \Ess \phi_i = \{1, \ldots, p\}
\]
for some $p < m$. Let $q = d-m+p$, and define the $q$-ary function $s'$ as
\[
s' = s(\phi'_1, \ldots, \phi'_m, x_{p+1}, \ldots, x_{q}),
\]
where $\phi'_i = \phi_i(x_1, \ldots, x_p, \constf{0}, \ldots, \constf{0})$. Then
\[
\begin{split}
\lefteqn{s'(x_1, \ldots, x_p, \phi_{m+1}, \ldots, \phi_d)} \\
& = s(\phi'_1, \ldots, \phi'_m, x_{p+1}, \ldots, x_q)(x_1, \ldots, x_p, \phi_{m+1}, \ldots, \phi_d) \\
& \begin{split}
{}= s(&\phi'_1(x_1, \ldots, x_p, \phi_{m+1}, \ldots, \phi_d), \ldots, \phi'_m(x_1, \ldots, x_p, \phi_{m+1}, \ldots, \phi_d), \\ &x_{p+1}(x_1, \ldots, x_p, \phi_{m+1}, \ldots, \phi_d), \ldots, x_q(x_1, \ldots, x_p, \phi_{m+1}, \ldots, \phi_d))
\end{split} \\
& = s(\phi_1, \ldots, \phi_d)
= f.
\end{split}
\]
This contradicts the minimality of the $\mathcal{C}$-decomposition $f = s(\phi_1, \ldots, \phi_d)$.
\end{proof}

A minimal $\mathcal{C}$-decomposition $(g, \phi_1, \ldots, \phi_d)$ of $f$ is called \emph{optimal,} if the cardinality $\card{\range (\phi_1, \ldots, \phi_d)}$ of the range of the inner functions is the smallest possible among all minimal $\mathcal{C}$-decompositions of $f$, and this smallest cardinality is called the \emph{$\mathcal{C}$-range degree} of $f$, denoted $\deg^\mathrm{r}_\mathcal{C} f$.

\begin{lemma}
\label{Crdegree}
If $f \subf[\mathcal{C}] g$ and $\deg_\mathcal{C} f = \deg_\mathcal{C} g$, then $\deg^\mathrm{r}_\mathcal{C} f \leq \deg^\mathrm{r}_\mathcal{C} g$.
\end{lemma}
\begin{proof}
Let $\deg_\mathcal{C} f = \deg_\mathcal{C} g = d$, $\deg^\mathrm{r}_\mathcal{C} g = r$, and let $(s, \phi_1, \ldots, \phi_d)$ be an optimal $\mathcal{C}$-decomposition of $g$. We have that $f = g(h_1, \ldots, h_n)$ for some $h_1, \ldots, h_n \in \mathcal{C}$, and so
\[
f = s(\phi_1, \ldots, \phi_d)(h_1, \ldots, h_n) = s(\phi'_1, \ldots, \phi'_d),
\]
where $\phi'_i = \phi_i(h_1, \ldots, h_n)$, and therefore $(s, \phi'_1, \ldots, \phi'_d)$ is a minimal $\mathcal{C}$-de\-com\-po\-si\-tion of $f$. Since $\range (\phi'_1, \ldots, \phi'_d) \subseteq \range (\phi_1, \ldots, \phi_d)$, we have that $\deg^\mathrm{r}_\mathcal{C} f \leq r$.
\end{proof}

\begin{corollary}
$\mathcal{C}$-equivalent functions have the same $\mathcal{C}$-range degree.
\end{corollary}

\begin{lemma}
\label{permif}
If $(g, \phi_1, \ldots, \phi_m)$ is an optimal $\mathcal{C}$-decomposition of $f$, then for every permutation $\sigma$ of $\{1,\ldots,m\}$, there is a function $g'$ such that $(g', \phi_{\sigma(1)},\linebreak[0] \ldots,\linebreak[0] \phi_{\sigma(m)})$ is an optimal $\mathcal{C}$-decomposition of $f$.
\end{lemma}
\begin{proof}
For any permutation $\sigma$ of $\{1, \ldots, d\}$,
\[
s(x_{\sigma^{-1}(1)}, \ldots, x_{\sigma^{-1}(d)})(\phi_{\sigma(1)}, \ldots, \phi_{\sigma(d)}) = s(\phi_1, \ldots, \phi_d).
\]
Thus, if $(s, \phi_1, \ldots, \phi_d)$ is an optimal $\mathcal{C}$-decomposition of $f$, then so is $(s(x_{\sigma^{-1}(1)},\linebreak[0] \ldots,\linebreak[0] x_{\sigma^{-1}(d)}), \phi_{\sigma(1)}, \ldots, \phi_{\sigma(d)})$.
\end{proof}

\subsection{Descending chain condition for $\subc[\mathcal{B}_1]$}

If $s = s(\psi_1, \ldots, \psi_n)$ for some essentially unary functions $\psi_1, \ldots, \psi_n$ such that the restriction of $(\psi_1, \ldots, \psi_n)$ to $S = \range \psi_1 \times \dots \times \range \psi_n = \range (\psi_1, \ldots, \psi_n)$ is the identity function on $S$, then we say that $s$ \emph{retracts} to $S$ and we call $(\psi_1, \ldots, \psi_n)$ a \emph{retraction map.}

\begin{lemma}
\label{psilemma}
Assume that $f = s(\phi_1, \ldots, \phi_d)$ is an optimal $\mathcal{B}_1$-decomposition. Then there is a function $s'$ such that $f = s'(\phi_1, \ldots, \phi_d)$ and $s'$ retracts to $\range \phi_1 \times \dots \times \range \phi_n$.
\end{lemma}
\begin{proof}
 For $i = 1, \ldots, d$, let $\psi_i$ be the essentially unary function defined as
\[
\psi_i(\vect{a}) =
\begin{cases}
a_i, & \text{if $a_i \in \range \phi_i$,} \\
\phi_i(\vect{0}), & \text{if $a_i \notin \range \phi_i$,}
\end{cases}
\]
and let $s' = s(\psi_1, \ldots, \psi_d)$. Then $f = s'(\phi_1, \ldots, \phi_d)$, $s' = s'(\psi_1, \ldots, \psi_d)$, $\range \psi_i = \range \phi_i$ for $1 \leq i \leq d$, and the restriction of $(\psi_1, \ldots, \psi_d)$ to $\range \phi_1 \times \dots \times \range \phi_d$ is the identity function.
\end{proof}

In the proof of the next proposition, we will make use of the following consequence of Rado's Theorem, also proved by Foldes and Lehtonen \cite{FL}.

\begin{theorem}
\label{matrix}
Let the columns of a $p \times q$ matrix $M$ over any field be partitioned into $n$ blocks, $M = [M_1, \ldots, M_n]$ ($p \leq n$). The following are equivalent.
\begin{enumerate}[(i)]
\item All $p \times p$ submatrices of $M$ with columns from distinct blocks $M_i$ are singular.
\item There is an invertible matrix $Q$ and an integer $m \geq 1$ such that in $QM = [QM_1, \ldots, QM_n]$, there are $m$ rows which are null in all but at most $m-1$ blocks $QM_i$.
\end{enumerate}
\end{theorem}

\begin{proposition}
\label{B1eqs}
If $f = s(\phi_1, \ldots, \phi_d)$ is an optimal $\mathcal{B}_1$-decomposition, then $\range (\phi_1,\linebreak[0] \ldots,\linebreak[0] \phi_d) = \range \phi_1 \times \dots \times \range \phi_d$. Furthermore, if $s$ retracts to $\range (\phi_1,\linebreak[0] \ldots,\linebreak[0] \phi_d)$, then $f \fequiv[\mathcal{B}_1] s$.
\end{proposition}
\begin{proof}
We call a function $f \in \mathcal{B}_1$ \emph{wide,} if it is not quasilinear. A wide function is essentially unary and its range contains at least three elements.

Let $f = s(\phi_1, \ldots, \phi_d)$ be an optimal $\mathcal{B}_1$-decomposition of an $n$-ary function $f$. By Lemma \ref{permif} we may assume that $\phi_1, \ldots, \phi_p$ are quasilinear and $\phi_{p+1}, \ldots, \phi_d$ are wide.
Denote $\Phi = \range \phi_1 \times \dots \times \range \phi_d$, $w = d-p$, $q = n-w = n-d+p$. Since for any permutation $\sigma$ of $\{1,\ldots,n\}$ and for any clone $\mathcal{C}$, $f \fequiv[\mathcal{C}] f(x_{\sigma(1)}, \ldots, x_{\sigma(n)})$, we may assume that, for $i = 1, \ldots, w$, $\Ess \phi_{p+i} = \{q+i\}$; let $\phi_{p+i} = \xi_{p+i}(x_{q+i})$ for a unary function $\xi_{p+i}$. We may assume that the quasilinear inner functions are of the form
\begin{equation}
h_1(x_1) \oplus \dots \oplus h_n(x_n), \quad \text{$h_j(0) = 0$ for $1 \leq j \leq n$.}
\label{canform}
\end{equation}
For, if $\phi_i$ has the standard form $\phi_i = g(h_1(x_1) \oplus \dots \oplus h_n(x_n))$, then
\[
\begin{split}
f &= s(\phi_1, \ldots, \phi_i, \ldots, \phi_d) \\
&= s(\phi_1, \ldots, \phi_{i-1}, g(h_1(x_1) \oplus \dots \oplus h_n(x_n)), \phi_{i+1}, \ldots, \phi_d) \\
&= s'(\phi_1, \ldots, \phi_{i-1}, h_1(x_1) \oplus \dots \oplus h_n(x_n), \phi_{i+1}, \ldots, \phi_d),
\end{split}
\]
where $s' = s(x_1, \ldots, x_{i-1}, \tilde{g}(x_i), x_{i+1}, \ldots, x_d)$ and $\tilde{g}$ is any permutation of $A$ whose restriction to $\{0,1\}$ coincides with $g$. Note that $s \fequiv[\mathcal{B}_1] s'$. Also, for the wide inner functions, we may assume that for $p + 1 \leq i \leq d$, $\range \phi_i = \nset{r_i}$, where $r_i = \card{\range \phi_i}$, and $\phi_i(\vect{0}) = 0$. We can make this condition hold with some suitable permutations $\sigma_i$ of $A$:
\[
\begin{split}
f &= s(\phi_1, \ldots, \phi_d) \\
&= s(x_1, \ldots, x_p, \sigma_{p+1}^{-1}(x_{p+1}), \ldots, \sigma_d^{-1}(x_d))(\phi_1, \ldots \phi_p, \sigma_{p+1} \circ \phi_{p+1}, \ldots, \sigma_d \circ \phi_d)
\end{split}
\]
and $s \fequiv[\mathcal{B}_1] s(x_1, \ldots, x_p, \sigma_{p+1}^{-1}(x_{p+1}), \ldots, \sigma_d^{-1}(x_d))$.

We denote $\zeta = \chi_{\{1\}}$.
If $s$ retracts to $\Phi$ then we also have that for $1 \leq i \leq p$,
\begin{equation}
s(x_1, \ldots, x_{i-1}, \zeta(x_i), x_{i+1}, \ldots, x_d) = s.
\label{eqzeta}
\end{equation}

We will now show that $\range (\phi_1, \ldots, \phi_d) = \Phi$.
For $1 \leq i \leq w$, let $S_{q+i}$ be a transversal of $\KER \xi_{p+i}$ such that $0 \in S_{q+i}$.
For $1 \leq i \leq p$, we may choose $\phi_i = h^i_1(x_1) \oplus \dots \oplus h^i_n(x_n)$ such that for $1 \leq j \leq w$, $h^i_{q+j}(a) = 0$ for every $a \in S_{q+j}$. For, assume that for some $1 \leq i \leq p$, $1 \leq j \leq w$, $a \in S_{q+j}$, we have $h^i_{q+j}(a) = 1$. Denote by $[a]$ the equivalence class of $a$ in $\KER \xi_{p+j}$, denote $\alpha = \xi_{p+j}(a)$, and define $\tilde{h}^i_{q+j} = h^i_{q+j} \oplus \chi_{[a]}$, i.e.,
\[
\tilde{h}^i_{q+j}(x) =
\begin{cases}
h^i_{q+j}(x), & \text{if $x \notin [a]$,} \\
h^i_{q+j}(x) \oplus 1, & \text{if $x \in [a]$.}
\end{cases}
\]
Then $\tilde{h}^i_{q+j}(a) = 0$. Now let
\[
\begin{split}
\phi'_i &= h^i_1(x_1) \oplus \dots \oplus h^i_{q+j-1}(x_{q+j-1}) \oplus \tilde{h}^i_{q+j}(x_{q+j}) \\
& \phantom{{}=h^i_1(x_1)} \oplus h^i_{q+j+1}(x_{q+j+1}) \oplus \dots \oplus h^i_n(x_n) \\
&= \phi_i \oplus \chi_{[a]}(x_{q+j}),
\end{split}
\]
and let $s' = s(x_1, \ldots, x_{i-1}, \rho, x_{i+1}, \ldots, x_d)$ with $\rho = \zeta(x_i) \oplus \chi_{\{\alpha\}}(x_{p+j})$. Then
\[
\begin{split}
\lefteqn{s'(\phi_1, \ldots, \phi_{i-1}, \phi'_i, \phi_{i+1}, \ldots, \phi_d)} \\
&= s(x_1, \ldots, x_{i-1}, \rho, x_{i+1}, \ldots, x_d)(\phi_1, \ldots, \phi_{i-1}, \phi'_i, \phi_{i+1}, \ldots, \phi_d) \\
&= s(\phi_1, \ldots, \phi_{i-1}, \rho(\phi_1, \ldots, \phi_{i-1}, \phi'_i, \phi_{i+1}, \ldots, \phi_d), \phi_{i+1}, \ldots, \phi_d) = f,
\end{split}
\]
because
\begin{multline*}
\rho(\phi_1, \ldots, \phi_{i-1}, \phi'_i, \phi_{i+1}, \ldots, \phi_d) \\
= \zeta(\phi'_i) \oplus \chi_{\{\alpha\}}(\phi_{p+j})
= \zeta(\phi_i \oplus \chi_{[a]}(x_{q+j})) \oplus \chi_{\{\alpha\}}(\xi_{p+j}(x_{q+j})) \\
= \phi_i \oplus \chi_{[a]}(x_{q+j}) \oplus \chi_{[a]}(x_{q+j})
= \phi_i.
\end{multline*}
If $s$ retracts to $\Phi$, then $s \fequiv[\mathcal{B}_1] s'$, because
\begin{multline*}
s'(x_1, \ldots, x_{i-1}, \rho, x_{i+1}, \ldots, x_d) \\
= s(x_1, \ldots, x_{i-1}, \rho(x_1, \ldots, x_{i-1}, \rho, x_{i+1}, \ldots, x_d), x_{i+1}, \ldots, x_d)
= s,
\end{multline*}
where the last equality holds by Equation \eqref{eqzeta} and since
\begin{multline*}
\rho(x_1, \ldots, x_{i-1}, \rho, x_{i+1}, \ldots, x_d)
= \zeta(\rho) \oplus \chi_{\{\alpha\}}(x_{p+j}) \\
= \zeta(\zeta(x_i) \oplus \chi_{\{\alpha\}}(x_{p+j})) \oplus \chi_{\{\alpha\}}(x_{p+j})
= \zeta(x_i).
\end{multline*}
Repeating this procedure, we will obtain a $\mathcal{B}_1$-decomposition $f = \tilde{s}(\tilde{\phi}_1,\linebreak[0] \ldots,\linebreak[0] \tilde{\phi}_p,\linebreak[0] \phi_{p+1},\linebreak[0] \ldots,\linebreak[0] \phi_d)$, where
\begin{equation}
h^i_{q+j}(a) = 0 \qquad \text{for every $1 \leq i \leq p$, $1 \leq j \leq w$, $a \in S_{q+j}$.}
\label{form3}
\end{equation}
In other words, the restrictions of $\tilde{\phi}_1, \ldots, \tilde{\phi}_p$ into $A^q \times S_{q+1} \times \dots \times S_n$ do not depend on the variables $x_{q+1}, \ldots, x_n$. Furthermore, the $\mathcal{B}_1$-decomposition $f = \tilde{s}(\tilde{\phi}_i, \ldots, \tilde{\phi}_p, \phi_{p+1}, \ldots, \phi_d)$ is optimal, and if $s$ retracts to $\Phi$ then $s \fequiv[\mathcal{B}_1] \tilde{s}$.

Thus, we can assume that the quasilinear inner functions have the form of Equation \eqref{form3}. We then consider the restrictions of the quasilinear inner functions $(\phi_1, \ldots, \phi_p)$ to $A^{q} \times S_{q+1} \times \dots \times S_n$. We may now assume that $\bigcup_{i=1}^p \Ess \phi_i = \{1, \ldots, q\}$.

We present a system of $p$ quasilinear functions $\phi_1, \ldots, \phi_p$ in the form of Equation \eqref{canform} ($\phi_i = h^i_1(x_1) \oplus \dots \oplus h^i_q(x_q)$) with $\bigcup_{i=1}^p \Ess \phi_i \subseteq \{1, \ldots, q\}$ as a matrix $M$ over the two-element field $\{0,1\}$ as follows. The rows of $M$ are indexed by $\{1,\ldots,p\}$, and the columns are indexed by $C = \{1, \ldots, q\} \times \{1, \ldots, k-1\}$. We let $M(i,(j,a)) = h^i_j(a)$. (Note that we are assuming that $h^i_j(0) = 0$ for all $i$ and $j$, so this information need not be encoded in $M$.) We then partition $C$ into $q$ blocks as $\Pi = \{C_1, \ldots, C_q\}$, where $C_j = \{j\} \times \{1, \ldots, k-1\}$. The elementary row operations (permutation of rows, addition of one row to another) correspond to permutation of the $\phi_i$'s and substitution of $\phi_i \oplus \phi_j$ for $\phi_i$ for some $i \neq j$.

The modulo $2$ sum of quasilinear functions of the form of Equation \eqref{canform} is again of this form: $\phi_i \oplus \phi_j = (h^i_1 \oplus h^j_1)(x_1) \oplus \dots \oplus (h^i_n \oplus h^j_n)(x_n)$ (note that $\chi_S \oplus \chi_{S'} = \chi_{S \bigtriangleup S'}$, where $\bigtriangleup$ denotes the symmetric difference). If $\phi_1, \ldots, \phi_m$ are functionally independent then so are also $\phi_1, \ldots, \phi_{i-1}, \phi_i \oplus \phi_j, \phi_{i+1}, \ldots, \phi_m$ for $i \neq j$.

Let $i \neq j$, and define the $d$-ary function $s'$ as
\[
s' = s(x_1, \ldots, x_{i-1}, \tau, x_{i+1}, \ldots, x_d),
\]
where $\tau = \zeta(x_i) \oplus \zeta(x_j)$.
Note that $\tau$ is a quasilinear function of the form of Equation \eqref{canform}. Then
\[
s'(\phi_1, \ldots, \phi_{i-1}, \phi_i \oplus \phi_j, \phi_{i+1}, \ldots, \phi_d)
= s(\phi_1, \ldots, \phi_d) = f.
\]
If $s$ retracts to $\Phi$, then we also have that $s = s'(x_1, \ldots, x_{i-1}, \tau, x_{i+1}, \ldots, x_d)$ by Equation \eqref{eqzeta}, so $s \fequiv[\mathcal{B}_1] s'$.

Since $\phi_1, \ldots, \phi_p$ are part of an optimal $\mathcal{B}_1$-decomposition $f = s(\phi_1, \ldots, \phi_d)$, it follows by Lemma \ref{Cdecev} that $M$ does not satisfy condition (ii) of Theorem \ref{matrix}, for otherwise we would have a $\mathcal{B}_1$-decomposition $f = s'(\phi'_1, \ldots, \phi'_d)$, where some $m$ inner functions only depend on $m-1$ variables, contradicting the optimality of the given $\mathcal{B}_1$-decomposition. Therefore Theorem \ref{matrix} implies that there is a set $D = \{(c_1, a_1), \ldots, (c_p, a_p)\}$ of $p$ columns of $M$ that is a partial transversal of $\Pi$ such that the square submatrix $R = M[\nset{p}, D]$ is nonsingular. Thus the range of $(\phi_1, \ldots, \phi_p)$ is the whole of $\{0,1\}^p$. Since the wide inner functions $\phi_{p+1}, \ldots, \phi_d$ do not depend on the first $q$ variables, we can now conclude that $\range (\phi_1, \ldots, \phi_d) = \Phi$.

We still have to show that if $s$ retracts to $\Phi$ then $s \subf[\mathcal{B}_1] f$. Assume that $s = s(\psi_1, \ldots, \psi_d)$, where $(\psi_1, \ldots, \psi_d)$ is a retraction map with range $\Phi$. For $i = 1, \ldots, q$, define the subset $S_i$ of $A$ as
\[
S_i =
\begin{cases}
\{0, a_j\}, & \text{if $i = c_j$,} \\
\{0\}, & \text{if $i \notin \{c_1, \ldots, c_p\}$.}
\end{cases}
\]
Recall that for $i = 1, \ldots, w$, $S_{q+i}$ was defined as a transversal of $\KER \xi_{p+i}$. It is not difficult to see that $S = S_1 \times \dots \times S_n$ is a transversal of $\KER (\phi_1, \ldots, \phi_d)$. Let the $p \times p$ matrix $Q = (q_{ij})$ be the inverse of $R$. Then the inverse mapping of $\phi|_S$ is $\delta = (\delta_1, \ldots, \delta_n)$, where
\[
\delta_i =
\begin{cases}
g_i(q_{j1}x_1 \oplus \dots \oplus q_{jp}x_p), & \text{if $i = c_j$,} \\
\constf{0}, & \text{if $i \in \{1, \ldots, q\} \setminus \{c_1, \ldots, c_p\}$,} \\
\xi^{-1}_{p+j}(x_{q+j}), & \text{if $i = q+j$ for some $1 \leq j \leq w$,}
\end{cases}
\]
where $g_i$ is the map $0 \mapsto 0$, $1 \mapsto a_j$. Letting $\delta'_i = \delta_i(\psi_1, \ldots, \psi_p)$ we have that $s = f(\delta'_1, \ldots, \delta'_n)$ and $\delta'_1, \ldots, \delta'_n \in \mathcal{B}_1$.
\end{proof}

\begin{proposition}
\label{B1equiv}
If $f \subf[\mathcal{B}_1] g$, $\deg_{\mathcal{B}_1} f = \deg_{\mathcal{B}_1} g$, and $\deg^\mathrm{r}_{\mathcal{B}_1} f = \deg^\mathrm{r}_{\mathcal{B}_1} g$, then $f \fequiv[\mathcal{B}_1] g$.
\end{proposition}
\begin{proof}
Let $g = s(\phi_1, \ldots, \phi_d)$ be an optimal $\mathcal{B}_1$-decomposition. We may assume that $s = s(\psi_1, \ldots, \psi_d)$ where $(\psi_1, \ldots, \psi_d)$ is a retraction map with $\range (\psi_1,\linebreak[0] \ldots,\linebreak[0] \psi_d) = \range (\phi_1, \ldots, \phi_d)$. Then by Proposition \ref{B1eqs}, $g \fequiv[\mathcal{B}_1] s$. We have that $f = g(h_1, \ldots, h_m)$ for some $h_1, \ldots, h_m \in \mathcal{B}_1$, and so $f = s(\phi_1, \ldots, \phi_d)(h_1, \ldots, h_m) = s(\phi'_1, \ldots, \phi'_d)$, where $\phi'_i = \phi_i(h_1, \ldots, h_m)$. This must be an optimal $\mathcal{B}_1$-de\-com\-po\-si\-tion of $f$ with $\range (\phi'_1, \ldots, \phi'_d) = \range (\phi_1, \ldots, \phi_d)$, and again by Proposition \ref{B1eqs}, $f \fequiv[\mathcal{B}_1] s$. By the transitivity of $\fequiv[\mathcal{B}_1]$, we have that $f \fequiv[\mathcal{B}_1] g$.
\end{proof}

\begin{theorem}
\label{B1chain}
$\subc[\mathcal{B}_1]$ satisfies the descending chain condition.
\end{theorem}
\begin{proof}
It follows from Lemma \ref{Cdegree}, Lemma \ref{Crdegree}, and Proposition \ref{B1equiv} that if $f \psubf[\mathcal{B}_1] g$, then either $\deg_{\mathcal{B}_1} f < \deg_{\mathcal{B}_1} g$ or $\deg_{\mathcal{B}_1} f = \deg_{\mathcal{B}_1} g$ and $\deg^\mathrm{r}_{\mathcal{B}_1} f < \deg^\mathrm{r}_{\mathcal{B}_1} g$. The $\mathcal{B}_1$-degree and the $\mathcal{B}_1$-range degree are nonnegative integers, and we cannot have an infinite descent in these parameters.
\end{proof}

\section{Infinite antichains in $\subc[\mathcal{B}_{k-2}]$}
\label{SecAC}

Assume that $\card{A} = k \geq 4$. For $n \geq 2$, define the $n$-ary function $f_n$ as
\[
f_n(a_1, \ldots, a_n) =
\begin{cases}
a_1, & \text{if $a_1 = \dots = a_n \neq k-1$,}\\
k-1, & \text{if $\card{\{i : a_i = k-1\}} = n - 1$,}\\
0, & \text{otherwise.}
\end{cases}
\]
It is clear that $\Ess f_n = \{1, \ldots, n\}$.

\begin{proposition}
\label{prop:Bk-2antichain}
For $n \neq m$, $f_n \incomp[\mathcal{B}_{k-2}] f_m$.
\end{proposition}
\begin{proof}
We say that a function $g \in \mathcal{B}_{k-2}$ is \emph{narrow} if $\card{\range g} \leq k-2$. We say that $g \in \mathcal{B}_{k-2}$ is \emph{wide} if $\card{\range g} > k-2$. Wide functions are essentially unary.

Assume that $n < m$. Suppose, on the contrary, that $f_m \subf[\mathcal{B}_{k-2}] f_n$. Then $f_m = f_n(\phi_1, \ldots, \phi_n)$ for some $\phi_1, \ldots, \phi_n \in \mathcal{B}_{k-2}$. Denote $\phi = (\phi_1, \ldots, \phi_n)$. We must have that for $1 \leq a \leq k-2$, $\phi(a, \ldots, a) = (a, \ldots, a)$, and so all the inner functions have a range of at least $k-2$ elements. The range of $\phi$ must also contain a vector with exactly $n-1$ elements equal to $k-1$. Thus, at least $n-1$ inner functions have a range of at least $k-1$ elements and are hence wide. Also, in order to obtain a subfunction of higher essential arity, at least one of the inner functions must be essentially at least binary and hence narrow. We conclude that $n-1$ inner functions are wide and one is narrow; by symmetry and without loss of generality, we may assume that $\phi_n$ is the narrow one with $\range{\phi_n} = \{1, 2, \ldots, k-2\}$. The other inner functions depend on one variable; assume that $\Ess \phi_i = \{\pi(i)\}$ for some $\pi : \{1, \ldots, n-1\} \rightarrow \{1, \ldots, m\}$. Consider the vector $\vect{v}$ with $\vect{v}(\pi(1)) = 1$ and $\vect{v}(j) = k-1$ for $j \neq \pi(1)$. We have that $\phi(\vect{v}) = (k-1, k-1, \ldots, k-1, x)$ for some $x \neq k-1$, so $\phi'_1(1) = k-1$, where $\phi'_1 = \phi_1(x_{\pi(1)}, \ldots, x_{\pi(1)})$. On the other hand, $\phi(1, \ldots, 1) = (1, \ldots, 1)$, so $\phi'_1(1) = 1$. We have reached a contradiction.

Suppose then, on the contrary, that $f_n \subf[\mathcal{B}_{k-2}] f_m$. Then $f_n = f_m(\phi_1,\linebreak[0] \ldots,\linebreak[0] \phi_m)$ for some $\phi_1, \ldots, \phi_m \in \mathcal{B}_{k-2}$. A similar argument as above shows that at least $m-1$ inner functions must be wide. Assume without loss of generality that the first $m-1$ inner functions are wide and $\Ess \phi_i = \{\pi(i)\}$ for some $\pi : \{1, \ldots, m-1\} \rightarrow \{1, \ldots, n\}$. The mapping $\pi$ must be injective. For, suppose on the contrary that $\pi(i) = \pi(j)$ for some $i \neq j$. Let $\vect{v} \in A^m$ with $\vect{v}(\pi(i)) = 1$ and $\vect{v}(h) = k-1$ for $h \neq \pi(i)$. Since $\phi(1, \ldots, 1) = (1, \ldots, 1)$, we have that $\phi(\vect{v})$ is a vector with at least two components equal to $1$, but this is a contradiction because $\phi(\vect{v})$ should be a vector with exactly $m-1$ components equal to $k-1$ and hence at most one component equal to $1$.

If $n < m-1$, then there is no such injective map $\pi$. If $n = m-1$, then $\pi$ is a permutation. Denote by $\vect{e}^n_i$ the $n$-vector whose $i$th component is $1$ and the other components are equal to $k-1$. Since for $1 \leq a \leq k-1$, $\phi(a, \ldots, a) = (a, \ldots, a)$, we must have that for $1 \leq i \leq m-1$, $\phi_i(\vect{v}) = 1$ for any $\vect{v}$ with $\vect{v}(\pi(i)) = 1$. Therefore we must have that for $1 \leq i \leq m-1$, $\phi(\vect{e}^n_{\pi(i)}) = \vect{e}^m_i$. But then also $\phi_m$ is wide; assume that $\Ess \phi_m = \{l\}$. But then $\phi(\vect{e}^n_l)$ would be a vector with two components equal to $1$, a contradiction.
\end{proof}

Unfortunately, the previous argument does not apply to $\mathcal{B}_1$ on a three-element base set. We have to treat this case differently. Assume that $\card{A} = k = 3$. For $n \geq 2$, define the $n$-ary function $g_n$ as
\[
g_n =
\begin{cases}
1, & \text{if $a_1 = \dots = a_n = 1$,} \\
2, & \text{if $\card{\{i : a_i = 2\}} = n-1$ and $\card{\{i : a_i = 0\}} = 1$,} \\
0, & \text{otherwise.}
\end{cases}
\]

\begin{proposition}
\label{prop:Bk-2antichaink3}
For $n \neq m$, $g_n \incomp[\mathcal{B}_1] g_m$. 
\end{proposition}
\begin{proof}
Let $n \neq m$ and suppose, on the contrary, that $g_m \subf[\mathcal{B}_1] g_n$. Then $g_m = g_n(\phi_1, \ldots, \phi_n)$ for some $\phi_1, \ldots, \phi_n \in \mathcal{B}_1$. Denote $\phi = (\phi_1, \ldots, \phi_n)$. Denote by $\vect{v}^n_i$ the $n$-vector whose $i$th component is equal to $0$ and the other components are equal to $2$. We have that $\phi(1,\ldots,1) = (1,\ldots,1)$ and $\phi(\vect{v}^m_i) = \vect{v}^n_{\pi(i)}$ for some $\pi : \{1,\ldots,m\} \to \{1,\ldots,n\}$.

The inner functions $\phi_i$ fall into two types: the quasilinear and the surjective. The surjective functions are essentially unary. We observe that if $\phi_i = \xi(x_j)$ for some unary surjective function $\xi$, then $\xi(1) = 1$ and so $\xi$ is either the identity function on $A$ or the mapping $0 \mapsto 2$, $1 \mapsto 1$, $2 \mapsto 0$. In other words, the surjective inner functions are either projections $x_j$ or \emph{negations} $\compl{x_j} = \xi'(x_j)$, where $\xi'$ is the latter of the two unary functions described above.

No projection occurs twice among the inner functions. For, if $\phi_r = \phi_s = x_t$ for some $t$ and $r \neq s$, then $\phi(\vect{v}^m_t)$ contains (at least) two $0$'s, which is not possible. Also, there is at most one negation among the inner functions. For, if $\phi_r = \compl{x_t}$, $\phi_s = \compl{x_u}$ for some $r \neq s$, then for any $l \notin \{t, u\}$, $\phi(\vect{v}^m_l)$ contains (at least) two $0$'s, which is not possible. If there are both a projection and a negation among the inner functions, then they depend on the same variable. For, if $\phi_r = x_t$, $\phi_s = \compl{x_u}$ for some $t \neq u$, then $\phi(\vect{v}^m_t)$ contains (at least) two $0$'s, which is again not possible. So, if there are surjections among the inner functions, then either they are all projections depending on distinct variables; or there are only one projection and one negation, which depend on the same variable; or there is only one negation and no projections.

If there is only one negation and no projections, say $\phi_r = \compl{x_t}$, then for any $l \neq t$, $\phi(\vect{v}^m_l) = \vect{v}^n_r$. Thus, for any $s \neq r$, $\range \phi_s = \{1,2\}$. Now, $\phi(\vect{v}^m_t) \neq \vect{v}^n_s$ for any $s$, because $\phi_r(\vect{v}^m_r) = 2$ and $0 \notin \range \phi_s$ for $s \neq r$. This is not possible.

If there are a projection and a negation, then assume without loss of generality that $\phi_1 = x_1$, $\phi_2 = \compl{x_1}$. Then for any $l$ and $3 \leq r \leq n$, $\phi_r(\vect{v}^m_l) = 2$, and so $\range \phi_r = \{1,2\}$. Let $\vect{u} = \phi(2,\ldots,2)$. Then $\vect{u} = (2, 0, u_3, \ldots, u_n)$ with $(2, \ldots, 2) \neq (u_3, \ldots, u_n) \in \{1,2\}^{n-2}$. We note that the value of $\phi_i$ changes from $1$ to $2$ (or vice versa) if the value of any variable changes from $1$ to $2$ (or vice versa) if and only if $i \in \{j : v_j = 1\}$. It is now easy to see that $\phi(0,0,0,2,2,\ldots,2) = \vect{v}^m_1$, a contradiction.

Thus, there is no negation among the inner functions. Assume then, without loss of generality, that $\phi_i = x_i$ for $1 \leq i \leq p < n$. Since $\phi(\vect{v}^m_1) = \vect{v}^n_1$, we have that $\range \phi_i = \{1,2\}$ for $p+1 \leq i \leq n$. But then $\phi(\vect{v}^m_{p+1}) \neq \vect{v}^n_j$ for every $j$, because $\phi_i(\vect{v}^m_{p+1}) = 2$ for $1 \leq i \leq p$ and $\phi_i(\vect{v}^m_{p+1}) \neq 0$ for $p+1 \leq i \leq n$. This is also an impossible situation.

Thus, either all inner functions are projections, or none of them is a projection. If $m < n$, then the number of distinct projections is less than the number of inner functions. If $m > n$, then there must be an essentially at least binary inner function in order to incorporate all $m$ essential variables of $g_m$. We have now established that none of the inner functions is surjective, essentially unary.

Consider now the only remaining case where all inner functions are quasilinear. Then there is a $j$ such that $\phi(\vect{v}^m_i) = \vect{v}^n_j$ for all $i$; say $j = 1$. Let $\vect{u} = \phi(2,\ldots,2)$, and let $V = \{i : \vect{u}(i) = 1\}$. We observe that the value of $\phi_i$ changes if the value of any variable changes from $1$ to $2$ (or vice versa) if and only if $i \in V$. It is now easy to see that $\phi(0,0,0,2,2,\ldots,2) = \vect{v}^n_1$, a contradiction.
\end{proof}

Propositions \ref{prop:Bk-2antichain} and \ref{prop:Bk-2antichaink3} can now be merged into one theorem that encompasses every $\card{A} = k \geq 3$.

\begin{theorem}
\label{Bk-2antichain}
For $\card{A} = k \geq 3$, there is an infinite antichain in $\subc[\mathcal{B}_{k-2}]$.
\end{theorem}

\end{document}